\documentclass[graybox]{svmult}

\usepackage{mathptmx}       
\usepackage{helvet}         
\usepackage{courier}        
\usepackage{type1cm}        
\usepackage{makeidx}         
\usepackage{graphicx}        
\usepackage{multicol}        
\usepackage[bottom]{footmisc}

\usepackage{amssymb,amsmath}

\providecommand{\abs}[1]{\lvert#1\rvert}
\newcommand{\exclude}[1]{}

\begin{document}

\title*{A De Bruijn--Erd\H{o}s theorem in graphs?}
\author{Va\v{s}ek Chv\'{a}tal}
\institute{Department of Computer Science and Software Engineering, Concordia University, Montr\'{e}al, Qu\'{e}bec, Canada, \email{chvatal@cse.concordia.ca}
}
\maketitle

\abstract{A set of $n$ points in the Euclidean plane determines at least $n$ distinct lines unless these $n$ points are collinear. In 2006, Chen and Chv\' atal asked whether the same statement holds true in general metric spaces, where the line determined by points $x$ and $y$ is defined as the set consisting of $x$, $y$, and all points $z$ such that one of the three points $x,y,z$ lies between the other two. The conjecture that it does hold true remains unresolved even in the special case where the metric space arises from a connected undirected graph with unit lengths assigned to edges. We trace its curriculum vitae and point out twenty-nine related open problems plus three additional conjectures.}
\section{Prehistory}
\label{sec:2}
It all started when Alain Guenoche and Bernard Fichet asked me if I wanted to come to their 
 third International Conference on Discrete Metric Spaces in Marseilles in September 1998. I like metric spaces and I love Marseilles, I replied, but I have no results I could present there. Never mind, they said magnanimously, come anyway. As I am not completely without shame, I then began racking my brain for something to talk about at the conference. A distant memory came to the rescue: As an undergraduate, I marvelled at the interpretation of families of sets as metric spaces provided by the Hamming metric on a family of indicator functions. Could a few combinatorial theorems be generalized to the realm of metric spaces? Dusting off my youthful ambition thirty years later, I circled around it till I settled on the project of looking for theorems of Euclidean geometry that might be generalized to arbitrary metric spaces.

\subsection{Lines and closure lines in metric spaces}
\label{subsec:lcl}
Saying that point $v$ in a Euclidean space lies between points $u$ and $w$ means that $v$ is an interior point of the line segment with endpoints $u$ and $w$; line $\overline{xy}$ is the set consisting of $x$, $y$, and all points $z$ such that one of the three points $x,y,z$ lies between the other two. These notions have straightforward extensions to arbitrary metric spaces:  In a space with metric $dist$, saying that point $v$ {\em lies between\/} points $u$ and $w$ means that $u,v,w$ are pairwise distinct and $dist(u,v)+dist(v,w)=dist(u,w)$; if {\em line\/} $L(xy)$ is defined as the set consisting of $x$, $y$, and all points $z$ such that one of the three points $x,y,z$ lies between the other two, then $L(xy)=\overline{xy}$ in the special case where $dist$ is the Euclidean metric.

This was the definition of lines in metric spaces that I hoped to use in extending a theorem or two of Euclidean geometry to arbitrary metric spaces. One candidate was the Sylvester--Gallai theorem~\cite{Syl,E82},
\begin{itemize}
\item[] Every non-collinear finite subset $V$ of the Euclidean plane such that $\abs{V} \ge 2$\\ includes two points such that the line determined by them passes through\\ no other point of $V$,
\end{itemize}
whose generalization would read
\begin{itemize}
\item[] In every finite metric space $(V,dist)$ such that $\abs{V} \ge 2$, some line  consists of only two points of $V$ or of all points of $V$.
\end{itemize}
This candidate flunked miserably: When $(V,dist)$ is the pentagon $C_5$ with the usual graph metric (in this case, $dist(x,y)=1$ when vertices $x,y$ are adjacent and $dist(x,y)=2$ when vertices $x,y$ are nonadjacent), $L(xy)$ consists of four vertices when $x,y$ are adjacent and it consists of three vertices when $x,y$ are nonadjacent.

Undaunted, I tried another tack: Let us define {\em closure line\/} $C(xy)$ as the smallest superset of $L(xy)$ such that $u,v\in C(xy)\Rightarrow L(uv)\subseteq C(xy)$. Just like lines $L(xy)$, closure lines $C(xy)$ are identical with Euclidean lines $\overline{xy}$ in the special case where the metric is Euclidean. Unlike lines $L(xy)$, closure lines $C(xy)$ did not flunk the Sylvester--Gallai test at once: I could not find a counterexample to the statement
\begin{itemize}
\item[] (SG) In every finite metric space $(V,dist)$ such that $\abs{V} \ge 2$, some closure line consists of only two points of $V$ or of all points of $V$.
\end{itemize} 
(In particular, $C_5$ is not a counterexample as each of its ten closure lines consists of all five vertices.)

\subsection{Sylvester--Gallai theorem in metric spaces?}
\label{subsec:sg}

Having formulated generalization (SG) of the Sylvester--Gallai theorem, I tried to prove it. 
The known proofs of the Sylvester--Gallai theorem~\cite{E43,C48,C61} did not help: I failed to adapt any of them to a proof of (SG). I considered the restricted version of (SG) where the metric spaces are induced by graphs: every connected undirected graph with vertex set $V$ induces the metric space $(V,dist)$ where $dist$ is the usual graph metric ($dist(u,v)$ standing for the number of edges in the shortest path from $u$ to $v$). This turned out to be easy: not only the restricted version of (SG), but even a stronger statement,
\begin{itemize}
\item[]  In every finite metric space $(V,dist)$ induced by a graph with at least two vertices, every closure line consists of only two points of $V$ or of all points of $V$,
\end{itemize} 
is valid. (The proof is a simple exercise: if $x$ and $y$ are adjacent twins, then $C(xy)=\{x,y\}$; else $C(xy)=V$.) To get more faith in the validity of (SG), I then tried to show that a counterexample could not be ridiculously small; plodding case analysis aided by computer search established that (SG) holds true for all metric spaces with at most nine points. Armed with this pathetic evidence, I presented the arrogant conjecture and related observations~\cite{Ch04} at the Marseilles meeting.

Over the next few years, I publicized the conjecture vigorously. I told it to anybody who would listen. I told it to first-class researchers and some of them may have taken a crack at it. I gave talks about it in different places. A mathematical luminary interrupted my lecture at Princeton to announce that he had a counterexample; a few minutes later he and the entire audience agreed that the example was not a counterexample. After the lecture; he proposed to me (now privately) a new counterexample; this, too, turned out to be false. Such episodes made me feel that the conjecture may have been not all that arrogant.

The conjecture remained unresolved till the fall of 2003.

\subsection{Enter Xiaomin}
\label{subsec:xm}

It was March 2000. There was a knock and when I opened my office door, there stood a young man who asked for a few minutes of my time. He explained to me his personal reasons for wanting to come to Rutgers as a graduate student in the middle of spring term and asked me if I could help by putting in a good word for him.

I said I sympathized, but as I didn't know him from Adam, I could not put in a good word for him. He replied that he anticipated this reaction and perhaps I could give him a test to get an idea of his mathematical abilities? As I was just about to leave for my graduate class in algorithms and data structures, I handed to him a copy of the midterm exam I was going to give in a few minutes and asked him to come back after class. He looked the exam over, asked for definitions of a couple of concepts he was unfamiliar with, and then we went our separate ways. When I returned and read his answers, my jaw dropped: they were a notch above those of the thirty students who had studied the material for half a term. It was only later and after much prodding from me that he reluctantly confessed to his high ranking in the Chinese Mathematical Olympiad. (China being a biggish country, I was much impressed, of course.)

I gave him a glowing recommendation, he was admitted, and the rest is history. His name was Xiaomin Chen.

\subsection{Sylvester--Gallai theorem in metric spaces!}
\label{subsec:2p}

In the fall of 2003, Xiaomin proved conjecture (SG).

A pivotal notion in Leroy Milton Kelly's celebrated short proof~\cite{C48}, \cite[Section 4.7]{C61}, \cite[Chapter 8]{G12} of the Sylvester--Gallai theorem is  the distance of a point from a line. This notion is unavailable in general metric spaces and yet echoes of Kelly's proof can be found in Chen's. Kelly minimizes the distance of point $b$ from line $\overline{ac}$ over all noncollinear triples $a,b,c$, which can be seen as choosing the flattest triangle with base $ac$ and apex $b$; Chen minimizes $dist(a,b) + dist(b,c) - dist(a,c)$, which can also be seen as choosing the flattest triangle with base $ac$ and apex $b$. Here, the following definitions are required to overcome a technical wrinkle: in a metric space: 
\begin{itemize}
\item{a {\em triangle\/} is a set of three points, none of which lies between the other two;}
\item{its three {\em edges\/} are its two-point subsets;}
\item{an edge is {\em simple\/} if no point lies between its two points;}
\item{a triangle is {\em simple\/} if all three of its edges are simple.}
\end{itemize}
Synopses of the two proofs are compared in Table~\ref{tab:1}. 
\begin{table}
\caption{Comparison of the two proofs}
\label{tab:1}
\resizebox{\textwidth}{!}{
\begin{tabular}{ll}
\hline\noalign{\smallskip}
Euclidean plane: Kelly & General metric space: Chen \\
\noalign{\smallskip}\svhline\noalign{\smallskip}
1. {\em If some three points of $V$ are noncollinear,} &1A. {\em If some three points are in no closure line,}\\
{\em then some line passes through only two points of $V$:}\phantom{xxx} &{\em then some simple triangle is in no closure line:}\
\\
 &\\
& if $a,b,c$ minimize \\
& $dist(a,b) + dist(b,c) + dist(a,c)$ \\
& over all triples of points in no closure line,\\
& then $\{a,b,c\}$ is a simple triangle.\\
& \\
&1B. {\em  If some simple triangle is in no closure line,}\\
&{\em  then some closure line consists of two points:}\\
\\
 if $a,b,c$ minimize  & if $a,b,c$ minimize \\
 the distance of point $b$ from line $\overline{ac}$, & $dist(a,b) + dist(b,c) - dist(a,c)$ \\
 over all noncollinear triples, & over all simple triangles,\\
 then $\overline{ac}$ passes through no third point of $V$.& then $C(ac)=\{a,c\}$.\\
& \\
2. {\em If every three points of $V$ are collinear,} &2. {\em If every three points are in some closure line,}\\
{\em then all points of $V$ are collinear.} &{\em then some closure line consists of all the points.}\\
\noalign{\smallskip}\hline\noalign{\smallskip}
\end{tabular}}
\end{table}

\noindent As for the devil in the details, the first part of Kelly's proof is crisp: if $\overline{ac}$ included three points $x,y,z$ with $y$  between $x$  and $z$, then the distance of $y$ from $\overline{bx}$ or the distance of $y$ from $\overline{bz}$ would be smaller than the distance of $b$ from $\overline{ac}$, a contradiction. By contrast, part 1B of Chen's proof is far from straightforward.

\section{A De Bruijn--Erd\H{o}s theorem in metric spaces?}
\label{sec: dbem}
When I was publicizing the conjecture that the Sylvester--Gallai theorem extends to metric spaces, Victor Klee and Benny Sudakov (independently of each other) suggested to me that other theorems on points and lines in the Euclidean plane might be eligible for a similar treatment. One of these is another well-known theorem,
\begin{theorem}\label{thm:predbe}
Every non-collinear finite subset $V$ of the Euclidean plane such that $\abs{V} \ge 2$ determines at least $\abs{V}$ distinct lines.
\end{theorem}
As Paul Erd\H os~\cite{E43} remarked in 1943, Theorem~\ref{thm:predbe} follows easily by induction from the Sylvester--Gallai theorem: 
\begin{quotation}
A line passing through only two points of $V$, point $x$ and another one, does not belong to the set of lines determined by $V-\{ x\}$. If this set includes at least $\abs{V}-1$ distinct lines, then $V$ determines at least $\abs{V}$ distinct lines; else, by the induction hypothesis, $V-\{ x\}$ is collinear, in which case lines $\overline{xy}$ with $y$ ranging over $V-\{ x\}$ are pairwise distinct.
\end{quotation}
In 2006, Xiaomin and I set out to investigate whether Theorem~\ref{thm:predbe} could be generalized to metric spaces. Candidate
\begin{itemize}
\item[] In every finite metric space $(V,dist)$ with at least two points, there are at least $\abs{V}$ distinct closure lines or some closure line  consists of all points of $V$
\end{itemize}
for such a generalization flunked badly:
\begin{theorem}\label{thm:cex}{\rm\cite[Theorem 7]{CheC}}
For every integer $n$ greater than $5$, there is a metric space on $n$
points where there are precisely $7$ distinct closure lines and each closure line consists of at most $n-2$ points.
\end{theorem}
Nevertheless, we could not find a counterexample with `closure lines' replaced by `lines'.
\begin{svgraybox}
\begin{conjecture}\label{con:dbem}
In every finite metric space with $n$ points such that $n\ge 2$, there are at least $n$ distinct lines or some line  consists of all $n$ points.
\end{conjecture}
\end{svgraybox}

\subsection{Terminology}
\label{subsec: term}
\subsubsection{Two De Bruijn--Erd\H{o}s theorems and one that is not}
\label{subsubsec: 2dbe}
Two joint results of  Nicolaas Govert de Bruijn and Paul Erd\H os share the name `De Bruijn--Erd\H{o}s theorem':
\begin{theorem}\label{thm:dbe1}{\rm\cite{DE48}}~
Let $m$ and $n$ be positive integers such that $m\ge 2$; let $V$ be
  a set of $n$ points; let $E$ be a family of $m$ subsets of $V$ such
  that every two distinct points of $V$ belong to precisely one member of $E$. Then
  $m\ge n$, with equality if and only if 
\begin{itemize}
\item $\mathcal L$ is of the type $\{p_1,p_2, \ldots , p_{n-1}\},\; \{p_1,p_n\},\; \{p_2,p_n\},\; \ldots ,\; \{p_{n-1},p_n\}$
\end{itemize}
or 
\begin{itemize}
\item $n=k(k-1)+1$ with each member of $\mathcal L$ containing $k$ points of $V$ and
each point of $V$ contained in $k$ members of $\mathcal L$.  
\end{itemize}
\end{theorem}
\begin{theorem}\label{thm:dbe2}{\rm\cite{DE51}}~
Let $k$ be a positive integer and let $G$ be an infinite graph. If every finite subgraph of $G$ is $k$-colourable, then $G$ is $k$-colourable.
\end{theorem}
Theorem~\ref{thm:predbe} is sometimes incorrectly referred to as the `De Bruijn--Erd\H{o}s theorem'. The confusion has no doubt originated from the fact that it is a special case of the far more powerful Theorem~\ref{thm:dbe1}. I apologize for having been one of the culprits perpetuating this error~\cite{CheC, ChiC} and I plead initial ignorance. 

So how should we refer to the proposition in Conjecture~\ref{con:dbem}? On the one hand, Theorem~\ref{thm:predbe} has no name; on the other hand, `De Bruijn--Erd\H{o}s theorem in metric spaces', while incorrect, is crisp and has been used for years. Let us stick to it.

\subsubsection{What is the meaning of `conjecture'?}
\label{subsubsec: meanconj}

For some people, conjecturing $X$ seems to imply commitment to the belief that $X$ is true. (My doctoral adviser Crispin St. John Alvah Nash-Williams would use the term `possible conjecture' when he was in doubt.) I am not one of these people: to me, conjecturing $X$ means that (i) I am interested in the truth value of $X$ and (ii) I have no counterexample. 
Still, when Xiaomin and I were writing our paper~\cite{CheC}, I phrased the problem as a question rather than a conjecture. I do not recall my reason for this cowardice. Anyway, what's in a name? That which we call a conjecture by any other name would \ldots ~\cite{Sha}.

\subsection{A logbook}
\label{subsec: log}

In the year 2000, the Government of Canada created a permanent program to establish 2000 research professorships called {\em Canada Research Chairs\/} in eligible degree-granting institutions across the country. I was privileged to hold first the Canada Research Chair in Combinatorial Optimization since my coming to Concordia in 2004 till 2011 and then the Canada Research Chair in Discrete Mathematics from 2011 till my retirement in 2014. Generous support from the program helped me in attracting brilliant postdocs and stellar visitors. 

We worked and played in the framework of a research group that met weekly for problem-solving sessions. To conform to the fashion of naming everything by an acronym, I dubbed it {\em ConCoCO\/} for Concordia Computational Combinatorial Optimization. ConCoCO became a meeting ground not only for my students and postdocs, but also for students and faculty from other Montreal universities and for short-time visitors. I was particularly touched by the steady an enthusiastic support of Luc Devroye, my long-time friend and former colleague at McGill.

How is all this relevant to the subject at hand? The range of topics discussed at ConCoCO extended far beyond Conjecture~\ref{con:dbem}, but the conjecture was its important part and the early results related to it came from ConCoCO.

\begin{itemize}

\item On 1 October 2006. Xiaomin and I first publicized our conjecture (oops! a question) along with a proof that every metric space on $n$  points has at least $\lg n$ distinct lines or a line consisting of all $n$ points.\\

\item On 2 October 2006, ConCoCO was inaugurated.\\ 

\item In September 2007, Ehsan Chiniforooshan joined ConCoCO as its first postdoc. He quickly became its cornerstone and then stayed on for two years. Two of our joint results (written up in May 2009) state that
\begin{itemize}
\item \cite[Corollary 1]{ChiC} in every metric space induced by a connected graph on $n$ vertices, there are $\Omega(n^{2/7})$ distinct lines or else some line consists of all $n$ vertices,
\item \cite[Theorem 3]{ChiC} in every metric space on $n$ points where each nonzero distance equals $1$ or $2$, there are $\Omega(n^{4/3})$ distinct lines and this bound is tight.
\end{itemize}
The lower bound $\Omega(n^{2/7})$ was later improved to $\Omega(n^{4/7})$: see Subsection~\ref{subsec:weak}. The lower bound $\Omega(n^{4/3})$ was later extended to metric spaces where each nonzero distance equals $1$, $2$, or $3$: see Subsection~\ref{subsec:metr}.\\

\item In August 2010, Google Scholar surprised me by telling me that we were not alone in the universe: others were interested in Conjecture~\ref{con:dbem}, too. I cannot resist the temptation to quote the opening of~\cite{JK} verbatim (except for the reference labels):
\begin{quotation}
In this paper, we present several results motivated by an open problem presented
by Chv\' atal in the problem session of IWOCA 2008. We study systems of lines
in metric spaces induced by graphs. Lines considered in this paper are sets of
vertices defined by a relation of betweenness, as introduced by Menger~\cite{M28}.
A line containing all the vertices is called a universal line. Similar properties,
concerning distances in graphs, are studied in metric graph theory, see a survey
by Bandelt and Chepoi~\cite{BC}.

The problem presented by Chv\' atal at IWOCA 2008 was originally conjectured by Klee and Wagon~\cite{KW}. It is a generalization of the De Bruijn--Erd\H os
Theorem~\cite{DE48}. The conjecture states the following:
\begin{itemize}
\item[] {\em Every graph with $n$ vertices defines at least $n$ different lines or it contains\\
a universal line.}
\end{itemize}
Klee and Wagon even stated this question about general discrete metric spaces,
but we consider only graphs. This problem is still open, see~\cite{CheC}.
\end{quotation}
The reader will draw her own conclusions.\\

\item In January 2011 I wrote to the legendary Maria Chudnovsky
\begin{quotation}
Some time during the period March -- September of this year, six of us are going to 
get together in Montreal for a week or so for a concentrated attack on a couple of 
problems concerning lines in finite metric spaces. The other five are strong 
mathematicians, but we have not made too much progress so far and I am beginning to 
wonder if the two problems are not beyond our reach.

And then it occurred to me that if only we could get you interested, you would end our 
misery by either solving the problems or, in the other case, certifying that they really 
are difficult.

So I wonder if there is any way to persuade you to join us. I would pay for your travel 
and hotel, of course, plus a little extra as a Concordia's visitor.
And we could schedule the workshop (not concurrent with but) adjacent to any event of 
your choice, such as the jazz festival or the film festival or anything else that you can 
find at\\
\phantom{are}{\tt http://www.tourisme-montreal.org/What-To-Do/Events}
\end{quotation}
To my great delight, she accepted this invitation.\\

Our meeting took place between June 3 and June 13. Eight of us were at its center: Laurent Beaudou (ConCoCO participant from April to November 2008 and then its postdoc from September 2009 to March 2010), Adrian Bondy (my long-time friend and favourite collaborator), Xiaomin, Ehsan, Maria, Nico Fraiman (ConCoCO participant since January 2011), Yori Zwols (Maria's multiple co-author and former doctoral student), and myself.\\

At our disposal we had the cozy  meeting room EV3.101 on the third floor of Concordia's Engineering, Computer Science, and Visual Arts Complex with its large windows overlooking rue St.~Catherine Street (as the signs in our bilingual Montreal used to say). There we would congregate every day from Saturday to Saturday at the crack of dawn and half an hour or so later, around 11AM, get down to business. In gruelling and greatly gratifying jam sessions we worked through the day nonstop, except for a short lunch break and a longer dinner break. (An important part of my job as the organizer was proposing new lunch and dinner venues every day.) These jam sessions went on till late at night and emails were flying from hotel room to hotel room well past midnight.\\

I brought two problems to the workshop. The first was to prove Conjecture~\ref{con:dbem} for metric spaces arising from connected chordal graphs and the second was to prove a generalization of the real De Bruijn--Erd\H os Theorem~\ref{thm:dbe1} in terms of $3$-uniform hypergraphs (see Subsection~\ref{subsec:3uh}). We solved the first problem as a warm-up in a day and half; the result was published much later as a note~\cite{BBC15} contending for the record of the least number of lines per author. The second problem turned out to be more difficult; it led us on an emotional roller coaster where exhilirating victories were rapidly turning into crushing defeats. By Friday evening we were sure of a consolation prize, a weaker and less elegant version of what I had proposed. Saturday morning we luxuriated in a proof of the whole thing. Saturday afternoon revealed a big hole in this proof, in the night from Saturday to Sunday the roller coaster rolled on, and on Sunday morning the hole was patched up. The workshop ended up on this fairy-tale note and the result was published some thirty months later~\cite{BBC13}.\\

\item Eleven weeks later, Yori joined ConCoCO as a postdoc for 2011--2012.\\

\item In January 2012, Cathryn Supko defected from McGill in order to become my M.Comp.Sc. student. She participated in ConCoCO with remarkable energy until her graduation in July 2014.\\
 
\item In April 2012, getting ready for the session "My Favorite Graph Theory Conjectures" of the June 2012 SIAM Conference on Discrete Mathematics in Halifax, I was reminded of a three years old irritant: Ehsan and I had  proved that $n$-point metric spaces where each nonzero distance equals $1$ or $2$ have $\Omega(n^{4/3})$ distinct lines, but we had not quite proved Conjecture~\ref{con:dbem} for these metric spaces. Even though our result provides a lower bound that is asymptotically far stronger than what the conjecture requires, it implies only that counterexamples to the conjecture, if any, include only finitely many metric spaces where each nonzero distance equals $1$ or $2$. I set out to remove this blemish. As it turned out, a variation on the arguments used in~\cite{ChiC} proved that the smallest counterexample to the restricted conjecture had to have at most $7$ points and then plodding case analysis took care of the rest. Ehsan was clearly entitled to a joint authorship, but he thought otherwise, and so I publicized it in {\tt arXiv}, and published it twenty-eight months later~\cite{Chv14}, as a single author.\\

\item In the same month, I found once again that interest in Conjecture~\ref{con:dbem} was not confined to our private group. This time, the outsiders' contribution was serious: Ida Kantor announced her talk at the forthcoming SIAM Conference on Discrete Mathematics. The talk presented results of her joint work with Bal\' azs Patk\' os on the conjecture restricted to the plane with the $\ell_1$ metric (Theorems~\ref{thm:kp1} and \ref{thm:kp2} in Subsection~\ref{subsec:metr}). Their admirable achievements made us feel less incestuous and their proof techniques inspired further work on the conjecture~\cite{ACH}.
\\

\item In the spring of 2013, I accepted two excellent postdocs, Pierre Aboulker and Rohan Kapadia, for the next academic year. (Later on, Rohan extended his stay by another year.) Around this time, Adrian told me of his planned visit to Montreal in late April and I jumped at the opportunity to try and re-create the magic atmosphere of July 2011.\\
\phantom{xxx} The new workshop took place on April 14 -- 27. Heraclitus was right, no man ever steps in the same river twice. Heavy teaching schedule prevented Maria from coming and Yori could not get away from his commitments, either. Ehsan, having kept in touch with us by email during the first week, drove to Montreal on Friday the 19th and back to Waterloo on Wednesday the 24th. Xiaomin was in constant touch with us by email from Shanghai. Pierre, the fresh PhD, joined us for the duration; Rohan, just before his own defense, could not. Cathryn, Laurent, Nico (and, of course, Adrian) were present all the time.\\
\phantom{xxx} Our aim was to improve the lower bound $\lg n$ on the number of distinct lines in $n$-vertex $3$-uniform hypergraphs where no line consists of all $n$ vertices~\cite[Theorem 4]{CheC}; see Subsection~\ref{subsec:3uh} for the definition. We kept improving the coefficient $1$ in front of the $\lg n$ little by little until Xiaomin's brilliant friend Peihan Miao, then a junior student in Shanghai Jiaotong University (and now a doctoral student at Berkeley) pushed it all the way up to $2-o(1)$. We published this result in~\cite{ABC14}.\\

\item In November 2013, Pierre and Rohan proved Conjecture~\ref{con:dbem} for metric spaces arising from connected distance-hereditary graphs~\cite{AK15}.\\

\item In the spring of 2014, Guangda Huzhang studied in his undergraduate thesis
at Shanghai Jiaotong University geometric dominant metric spaces and graphs.
His work was later expanded into a paper written jointly with Xiaomin, Peihan, and another of their 
friends, Kuan Yang~\cite{CHMY}. Some of their results are quoted here in Subsection~\ref{subsec:lines}.\\

\item ConCoCO held its last meeting on Thursday, 26 June 2014.\\

\item With my retirement on August 31, 2014, research into Conjecture~\ref{con:dbem} gained a new momentum:
\begin{itemize}
\item In the final four months of 2014, Pierre, Xiaomin, Guangda, Rohan, and Cathryn completed a project they had been working on since the beginning of the year. They published its results in~\cite{ACH}. Some of them are quoted here as Theorem~\ref{thm:glb} in Subsection~\ref{subsec:weak}, Theorem~\ref{thm:mlb} in Subsection~\ref{subsec:metr}, and Theorem~\ref{thm:psd} in Subsection~\ref{subsec:pseudo}.\\

\item On 28 January 2015, Pierre and his friends Guillaume Lagarde, David Malec, Abhishek Methuku, and Casey Tompkins posted on {\tt arXiv} their manuscript that was later published as~\cite{ALM}. Their results involve an analog of Conjecture~\ref{con:dbem} for a class of $3$-uniform hypergraphs, which is quoted here as Theorem~\ref{thm:alm} in Subsection~\ref{subsec:3uh}.\\

\item On 20 June, 2016, Pierre and his friends Martin Matamala, Paul Rochet, and Jos\' e Zamora posted on {\tt arXiv} their manuscript that was later published as~\cite{AMRZ}. They proved Conjecture~\ref{con:dbem} for metric spaces arising from a class of graphs that contains all connected chordal graphs and all connected distance-hereditary graphs. In addition, they proposed an intriguing variation on Conjecture~\ref{con:dbem}. Some of their results are quoted here in Subsection~\ref{subsec:inter}.\\

\end{itemize}

\end{itemize}

\section{A De Bruijn--Erd\H{o}s theorem in graphs?}
\label{sec: dbeg}
Every connected undirected graph $G$ gives rise to the metric space $M(G)$ by defining the distance between two vertices as the smallest number of edges in a path joining them. Let us not play at being pedants: Rather than talking of lines in $M(G)$, let us talk of lines in $G$. Conjecture~\ref{con:dbem} remains open even in the special case where the metric space arises from a graph (and isn't this fortunate, since otherwise how could I submit this piece to the collection entitled `{\em Graph Theory\/} Favorite Conjectures and Open Problems'?).

\begin{svgraybox}
\begin{conjecture}\label{con:dbeg}
In every finite connected graph with $n$ vertices such that $n\ge 2$, there are at least $n$ distinct lines or some line  consists of all $n$ vertices.
\end{conjecture}
\end{svgraybox}

\subsection{Terminology and notation}
\label{subsec: gterm}
All our graphs as well as metric spaces and related objects are finite (unless specified otherwise), and so we will skip the qualifier `finite' throughout the text. All our graphs are also undirected and connected (unless specified otherwise), and so we will skip these two qualifiers as well. To avoid the one-vertex graph (which has no lines at all), let us also agree that all our graphs have at least two vertices. 

We  let $\abs{G}$ denote the number of vertices in a graph $G$. A line in a graph $G$ is said to be {\em universal\/} if it consists of all $\abs{G}$ vertices. A graph $G$ is said to have {\em the DBE property\/} if it has at least $\abs{G}$ distinct lines or a universal line. In these terms and under our assumptions, Conjecture~\ref{con:dbeg} asserts that all graphs have the DBE property. 

Sometimes we write simply $uv$ for the unordered pair $\{u,v\}$ of distinct elements $u$ and $v$.

\subsection{A weaker lower bound attained by all graphs}
Proving that almost all graphs have $\Theta(n^2)$ distinct lines, whether they have a universal line or not, is an easy exercise. This is far more than Conjecture~\ref{con:dbeg} requires. When it comes to all graphs, we have only far less than Conjecture~\ref{con:dbeg} requires. 

\label{subsec:weak}
\begin{theorem}\label{thm:glb}{\rm \cite[Theorem 7.6]{ACH}}~
Every graph $G$ has $\Omega(\abs{G}^{4/7})$ distinct lines or a universal line.
\end{theorem}

\subsection{Special cases where the lower bound is attained}
\label{subsec:spec}
One way of making progress toward the proof of Conjecture~\ref{con:dbeg} is finding larger and larger classes of graphs with the DBE property. By now, we know three such classes:
\begin{itemize}
\item
\begin{theorem}\label{thm:bip} 
All bipartite graphs have the DBE property.
\end{theorem}
\vspace{0.1cm}
This theorem is just a simple observation: in a bipartite graph, $L(uv)$ is universal whenever $u$ and $v$ are adjacent. (For every vertex $w$, we have $\abs{\,dist(w,u)-dist(w,v)\,}\le 1$ and, since the graph is bipartite, $dist(w,u)\ne dist(w,v)$.) For its strengthening, see Theorem~\ref{thm:mzbip} in Subsection~\ref{subsec:inter}.
\vspace{0.2cm}
\item 
\begin{theorem}\label{thm:d2}{\rm\cite[special case of Theorem 1]{Chv14}}~
All graphs of diameter $2$ have the DBE property.
\end{theorem}
\vspace{0.2cm}
\item 
\begin{theorem}\label{thm:amrz}{\rm\cite[corollary of Theorem 2.1]{AMRZ}}~
All graphs that can be constructed from chordal graphs 
by repeated substitutions and gluing vertices  have the DBE property.
\end{theorem}
\vspace{0.1cm}
Theorem~\ref{thm:amrz} provides a common generalization of two previous results,
\begin{itemize}
\item {\em all chordal graphs have the DBE property\/}~\cite[Theorem 1]{BBC15} and 
\item {\em all distance-hereditary graphs have the DBE property\/}~\cite[Theorem 1]{AK15}.
\end{itemize}
For its strengthening, see Theorem~\ref{thm:amrz+} in Subsection~\ref{subsec:inter}.
\end{itemize}

\noindent Here are three challenges motivated by these three theorems:

A graph is called {\em bisplit\/} \cite{BHL} if its vertex set can be partitioned into stable sets $X$, $Y$, and $Z$ so that $Y\cup Z$ induces a complete bipartite graph. (Bipartite graphs are bisplit graphs with $Z=\emptyset$.)
\begin{svgraybox}
\begin{problem}\label{prb:bs} 
Prove that all bisplit graphs have the DBE property.
\end{problem}
\end{svgraybox}
\begin{svgraybox}
\begin{problem}\label{prb:d3} 
Prove that all graphs of diameter $3$ have the DBE property.
\end{problem}
\end{svgraybox}
\noindent It is known that every  graph $G$ of diameter $3$ has at least $\abs{G}/15$ distinct lines or a universal line: more generally, every  graph $G$ of diameter $k$ has at least $\abs{G}/5k$ distinct lines or a universal line (see Theorem~\ref{thm:diam+} in Subsection~\ref{subsec:metr}.)\\

The {\em house\/} is the complement of the chordless path on five vertices; 
a {\em hole\/} is a chordless cycle with at least five vertices; 
the {\em domino\/} is the cycle on six vertices with one long and no short chord. An 
{\em HHD-free graph\/}~\cite{HK} contains no  house, no hole, and no domino as an induced subgraph. 
All graphs featured in Theorem~\ref{thm:amrz} are HHD-free, but not all HHD-free graphs can be constructed as in Theorem~\ref{thm:amrz}. For example, start with the $C_4$ that has vertices $ v_1, v_2, v_3, v_4$ and edges $v_1v_2, v_2v_3, v_3v_4, v_4v_1$. Then, for each of the two $i=1,2$, substitute a clique $\{a_i,c_i,e_i\}$ for $v_i$, add vertices $b_i$, $d_i$, and add edges $a_ib_i$, $b_ic_i$, $c_id_i$, $d_ie_i$.
\begin{svgraybox}
\begin{problem}\label{prb:hhd} 
Prove that all HHD-free graphs have the DBE property.
\end{problem}
\end{svgraybox}
Theorem~\ref{thm:amrz} would follow from Theorem 1 of \cite{BBC15} if it were known that substitution preserves the DBE property and that gluing vertices preserves the DBE property. As for the former proposition, it is not known that substitution preserves the DBE property even in the special case where the (not necessarily connected) graph that is being substituted for a vertex has only two vertices.
\begin{svgraybox}
\begin{problem}\label{prb:adj} 
Prove that splitting a vertex into adjacent twins preserves the DBE property.
\end{problem}
\end{svgraybox}
\begin{svgraybox}
\begin{problem}\label{prb:non} 
Prove that splitting a vertex into nonadjacent twins preserves the DBE property.
\end{problem}
\end{svgraybox}

\begin{svgraybox}
\begin{problem}\label{prb:dglue} 
Prove that gluing vertices preserves the DBE property.
\end{problem}
\end{svgraybox}
\noindent Ehsan Chiniforooshan and Xiaomin Chen [personal communication] solved a special case of Problem~\ref{prb:dglue}:
\begin{theorem}\label{thm:ehxi}
All graphs that can be constructed from graphs of diameter $2$ 
by repeatedly gluing vertices  have the DBE property.
\end{theorem}

\medskip

{\em Gallai graphs\/} (also known as $i$-triangulated graphs) are a common generalization of chordal graphs and bipartite graphs: every odd cycle of length 
at least five has at least two non-crossing chords. 
\begin{svgraybox}
\begin{problem}\label{prb:gal} 
Prove that all Gallai graphs have the DBE property.
\end{problem}
\end{svgraybox}
\noindent Since every Gallai graph with no clique-cutset is either a complete multipartite graph or else the join of a connected bipartite graph and a clique~\cite{Gal},
 Problem~\ref{prb:gal} is related to proving that
\begin{itemize}
\item[] (?) all graphs that can be constructed from graphs of diameter $2$ 
by repeated gluing along cliques have the DBE property,
\end{itemize}
which would strengthen Theorem~\ref{thm:ehxi}.\\

Theorems~\ref{thm:amrz} and~\ref{thm:ehxi} highlight the theme of building classes of graphs with the DBE property from prescribed classes by prescribed operations. One of the many additional variations on this theme goes as follows:
\begin{svgraybox}
\begin{problem}\label{prb:bip+} 
Prove that all graphs that can be constructed from bipartite graphs by repeated splitting of vertices into adjacent twins have the DBE property.
\end{problem}
\end{svgraybox}

\medskip
The family of {\em perfect graphs\/}~\cite{CRST} is, by definition, closed under taking induced subgraphs. Bisplit graphs, HHD-free graphs, and Gallai graphs are subfamilies of this family and they are also closed under taking induced subgraphs. This property seems irrelevant to graph metric; bisplit graphs, HHD-free graphs, and Gallai graphs are featured here just because they have been studied elsewhere and their structure is well understood. The last problem in this subsection concerns a possible strengthening of Conjecture~\ref{con:dbeg} for graphs in another family closed under taking induced subgraphs.

Every bridge in a graph defines a universal line, but not every universal line is defined by a bridge: for instance, the universal line in the wheel with five vertices is defined only by pairs of nonadjacent vertices. This graph and many other examples of bridgeless graphs with universal lines contain an induced subgraph isomorphic to $C_4$. Yori Zwols [personal communication] conjectured that the answer to the following question is `true':
\begin{svgraybox}
\begin{problem}\label{prb:yori} 
True or false? Every $C_4$-free graph $G$ has at least $\abs{G}$ distinct lines or a bridge.
\end{problem}
\end{svgraybox}
\noindent In March 2018, Martin Matamala and Jos\' e Zamora~\cite{MZ} proved his conjecture for bipartite graphs (see Theorem~\ref{thm:mzbip} in Subsection~\ref{subsec:inter}). 

\subsection{A red herring?}
\label{subsec:asym}

The lower bound $n$ in Conjecture~\ref{con:dbeg} is inherited from the more general Conjecture~\ref{con:dbem}. In the more general context of metric spaces, this bound (if at all valid) is tight (consider $n-1$ collinear 
points in the Euclidean plane and a point off their line). In the more restricted context of graphs, it may be so far from being tight as to be downright misleading:
\begin{svgraybox}
\begin{conjecture}\label{con:asym}
All graphs $G$ without a universal line have 
\begin{center}
$\Omega(\abs{G}^{4/3})$\ distinct lines.
\end{center}

\end{conjecture}
\end{svgraybox}
\noindent This conjecture emerged in our ConCoCO discussions of a theorem implying that 
all graphs $G$ of diameter at most $2$ have $\Omega(\abs{G}^{4/3})$ distinct lines~\cite[Theorem 3]{ChiC}.
As noted in~\cite[Theorem 3]{ChiC}, its lower bound is best possible: complete multipartite graphs with $\Theta(\abs{G}^{2/3})$ parts of sizes in $\Theta(\abs{G}^{1/3})$ have $\Theta(\abs{G}^{4/3})$ distinct lines and no universal line. There are many other graphs with these properties: Example 7.8 in ~\cite{ACH} exhibits arbitrarily large graphs $G$ with $\Theta(\abs{G}^{4/3})$ distinct lines, no universal line, and unbounded diameter. A neat variation on this theme has been pointed out by Xiaomin Chen: each graph in his class consists of $\Theta(n^{2/3})$ chordless cycles of lengths in $\Theta(n^{1/3})$ that, apart from a vertex common to all of them, are pairwise vertex-disjoint.\\

Conjecture~\ref{con:asym} is known to be valid for graphs of bounded diameter: here is a more general result.
\begin{theorem}\label{thm:asym3}{\rm \cite[Theorem 7.4]{ACH}}~
If $d(n)=o(n)$, then all graphs with $n$ vertices and diameter $d(n)$ have 
$\Omega((n/d(n))^{4/3})$ distinct lines.
\end{theorem}
Conjecture~\ref{con:asym} is also known to be valid for graphs where no line contains another (see Theorem~\ref{thm:gd} in Subsection~\ref{subsec:lines}). Its validity would not imply validity of Conjecture~\ref{con:dbeg}: it would imply only that counterexamples to Conjecture~\ref{con:dbeg}, if any, are finitely many. A plausible common strengthening of both conjectures goes as follows:
\begin{svgraybox}
\begin{conjecture}\label{con:nonasym}
For every graph $G$ with no universal line there is a complete multipartite graph with no universal line, as many vertices as $G$, and at most as many distinct lines as $G$.
\end{conjecture}
\end{svgraybox}
\noindent To see that Conjecture~\ref{con:nonasym} is indeed a common strengthening of both Conjecture~\ref{con:dbeg} and Conjecture~\ref{con:asym}, let $f(n)$ denote the smallest number of distinct lines in a complete multipartite graph with $n$ vertices and no universal line. Consider a complete $k$-partite graph $H$ with $n_i$ vertices in the $i$-th part and $n$ vertices altogether. If $H$ has no universal line, then $k\ge 3$ and $n_i\ne 2$ for all $i$, in which case $H$ has $\binom{k}{2}+\sum_{i=1}^k\binom{n_i}{2}$ distinct lines. To see that $f(n)\ge n$, assume without loss of generality that $n_i=1$ when $1\le i\le m$ and $n_i\ge 3$ when $m<i\le k$ for some $m$; then observe that
\[
\textstyle{\binom{k}{2}+\sum_{i=1}^k\binom{n_i}{2} \;=\; \binom{k}{2}+\sum_{i=m+1}^k\binom{n_i}{2} \;\ge\; k+\sum_{i=m+1}^k n_i \;=\; k+(n-m)\;\ge\; n.
}
\]
To see that $f(n)=\Omega(n^{4/3})$, observe that 

$$\textstyle{\binom{k}{2}+\sum_{i=1}^k\binom{n_i}{2}\ge \binom{k}{2}+k\binom{n/k}{2}\ge \frac{1}{2}\left(k^2+\frac{n^2}{k}\right)-n \ge 
\left(\frac{27}{32}\right)^{1/3}n^{4/3}-n.
}$$ 
\noindent Conjecture~\ref{con:nonasym} was suggested by recent experimental results of Yori Zwols. He computed the smallest number of distinct lines in graphs with at most $11$ vertices and no universal line. With a single exception, graphs attaining the minimum turned out to be complete multipartite. The exception, which is the complement of the Petersen graph, has $15$ distinct lines, just like the complete multipartite graphs $K_{3,3,4}$ and $K_{1,3,3,3}$.

\subsection{Families of lines}\label{subsec:lines}

Families of lines in graphs  have properties that
may seem outlandish to a visitor from a Euclidean space: for instance, the star with vertices $1,2,3,4$ and edges $12$, $13$, $14$ has lines $\{1,2,3\}$, $\{1,2,4\}$, $\{1,3,4\}$, $\{1,2,3,4\}$. Every two  lines in a Euclidean space  share at most one point, which is not the case in this example.
\begin{svgraybox}
\begin{problem}\label{prb:LG}
How difficult is it to recognize hypergraphs whose hyperedge set is the family of lines in some graph?
\end{problem}
\end{svgraybox}

A graph is {\em geometric dominant\/}~\cite{CHMY} if none of its lines contains another. In particular, if 
every two  lines in a graph share at most one vertex, then this graph is geometric dominant.
\begin{theorem}\label{thm:sgd}{\rm \cite[Theorem 4]{CHMY}}~
Every two  lines in a graph share at most one vertex if and only if this graph is complete or a path or $C_4$.
\end{theorem}
Geometric dominant graphs not specified in Theorem~\ref{thm:sgd} are called {\em nontrivial.\/} Nontrivial geometric dominant graphs with $n$ vertices may be hard to find when $n$ is small (the smallest one is the wheel with six vertices), but they are abundant when $n$ is large:
\begin{theorem}\label{thm:rand}{\rm \cite[Theorem 5]{CHMY}}~
If $p(n)^3n/\log n\rightarrow\infty$ and $(1-p(n))^2n/\log n\rightarrow\infty$
as $n\rightarrow\infty$, then the random graph ${\cal G}_{n,p(n)}$ is almost surely geometric dominant.
\end{theorem}
\begin{svgraybox}
\begin{problem}\label{prb:gd} 
Prove that all geometric dominant graphs have the DBE property.
\end{problem}
\end{svgraybox}
\begin{theorem}\label{thm:gd}{\rm \cite[Theorem 8]{CHMY}}~
All nontrivial geometric dominant graphs $G$ have $\Omega(\abs{G}^{4/3})$ distinct lines.
\end{theorem}
Proving that all nontrivial geometric dominant graphs have bounded diameter would make Theorem~\ref{thm:gd} a corollary of Theorem~\ref{thm:asym3}. Even a stronger statement might be true:
\begin{svgraybox}
\begin{problem}\label{prb:gd2}{\rm \cite[Question 1]{CHMY}}~ 
True or false? All nontrivial geometric dominant graphs have diameter $2$.
\end{problem}
\end{svgraybox}

\subsection{Equivalence relations}\label{subsec:equiv}

{\em Ceterum autem censeo Carthaginem delendam esse\/} (besides, I also believe that Carthage must be destroyed) was Cato the Elder's stock conclusion to all his speeches in the Roman Senate, irrespective of their topic. With similar persistence, Adrian Bondy liked to point out again and again in our ConCoCO discussions that our progress would get a great boost if we understood which equivalence relations $\equiv$ on the edge sets of $K_n$ arise from  graphs with $n$ vertices (or, more generally, from metric spaces on $n$ points) in the sense that $ab\equiv xy \Leftrightarrow L(ab)=L(xy)$.

Section 6 of~\cite{ACH} contains results on distinct pairs of vertices that define the same line. In its notation (Definition~6.2),
\begin{eqnarray*}
I(a,b) &=& \{z:\; \text{$z$ lies between $a$ and $b$}\},\\
O(a,b) &=& \{z:\; \text{$a$ lies between $z$ and $b$ or $b$ lies between $a$ and $z$}\}
\end{eqnarray*}
(so that $L(ab)=\{a,b\} \cup I(a,b)\cup O(a,b)$; in its terminology (Definitions~6.3 -- 6.5 and Lemma~6.9), a {\em parallelogram\/} is an ordered $4$-tuple $(a,b,c,d)$ of distinct vertices such that\\
\phantom{xxx}$\bullet\;$ $dist(a,b)=dist(c,d)$,\\
\phantom{xxx}$\bullet\;$ $dist(b,c)=dist(d,a)$,\\
\phantom{xxx}$\bullet\;$ $dist(a,c)=dist(b,d)=dist(a,b)+dist(b,c)$.
\begin{theorem}\label{thm:equiv}{\rm \cite[Lemma 6.6]{ACH}}~
Let $G$ be a   graph and let $e,f$ be distinct edges of the complete graph on the vertex set of $G$. If $L_G(e)=L_G(f)$, then the endpoints of $e$ can be labeled $a,b$ and the endpoints of $f$ can be labeled $c,d$ (possibly $b=c$) in such a way that\\ 
\phantom{} ($\alpha$)\, $b$ lies between $a$ and $c$;\ $c$ lies between $b$ and $d$; both $b,c$ lie between $a$ and $d$ or\\
\phantom{} ($\beta$) $(a,b,c,d)$ is a parallelogram and $I(a,b)=I(c,d)=\emptyset$ or\\
\phantom{} ($\gamma\;$) $(a,c,b,d)$ is a parallelogram and $O(a,b)=O(c,d)=\emptyset$.
\end{theorem}

\begin{svgraybox}
\begin{problem}\label{prb:geq}
How difficult is it to recognize equivalence relations $\equiv$ such that 
$ab\equiv xy$ if and only if $L_G(ab)=L_G(xy)$ for some graph $G$?
\end{problem}
\end{svgraybox}\noindent
Some of the candidates $\equiv$ partitioning edge sets of $K_n$ into classes $C_1$,\ldots $C_m$ are rejected by the following procedure.
\begin{tabbing}
{\bf Algorithm G:}\\
{\bf for}\, $i=1$ to $m$\, {\bf do}\, $L_i=$ the set of all endpoints of edges in $C_i$\, {\bf end}\\
{\bf while}\, \=there are pairwise distinct vertices $u,v,w$ and (not necessarily distinct)\\
\> subscripts $i,j$ such that $uv\in C_i$, $w\not\in L_i$, $vw\in C_j$, $u\in L_j$\\ 
{\bf do}\>add $w$ to $L_i$;\\ 
{\bf end}\\
{\bf if}\hspace{0.8cm}\=there are distinct $i,j$ such that $\abs{L_i}=\abs{L_j}=n$\\
{\bf then}\>{\bf return} message {\sc Does not arise from any  graph};\\ 
{\bf else}\>{\bf return} message {\sc Don't know};\\ 
{\bf end}
\end{tabbing}
For instance, given classes\ $C_1=\{12,23,34\}$\ and\ $C_2=\{13,24,14\},$ Algorithm G constructs $L_1=L_2=\{1,2,3,4\}$, and so it returns message {\sc Does not arise from any  graph}. Nevertheless, Algorithm G does not eliminate all inputs that do not arise from any  graph. For instance, given the partition into classes $C_1=\{14\}$, $C_2=\{24\}$, $C_3=\{34\}$, $C_4=\{12,23,13\}$ that does not arise from any  graph, Algorithm G constructs $L_1=\{1,4\}$, $L_2=\{2,4\}$, $L_3=\{3,4\}$, $L_4=\{1,2,3\}$, and so it returns message {\sc Don't know}. 

Correctness of Algorithm G follows from the observation that, for all graphs $G$ such that $L_G(uv)=L_G(xy)$ if and only if $uv$ and $xy$ belong to the same $C_i$, its {\bf while} loop maintains the invariant $uv\in C_i \;\Rightarrow\; L_G(uv)\supseteq L_i$.

\subsection{An interpolation}\label{subsec:inter}

Let us call an unordered pair $uv$ of vertices a {\em mighty pair} if $L(uv)$ is universal, let $\lambda(G)$ stand for the number of distinct lines in $G$, and let $\mu(G)$ stand for the number of mighty pairs in $G$. In this 
notation, Conjecture~\ref{con:dbeg}  states that 
$$\lambda(G)\ge \abs{G}\;\;\vee\;\; \mu(G)>0.$$ 
A stronger conjecture interpolates between the two operands of the disjunction:
\begin{svgraybox}
\begin{conjecture}\label{con:amrz}{\rm \cite[Conjecture 2.3]{AMRZ}}~
All graphs $G$ satisfy $\lambda(G)+\mu(G)\ge \abs{G}$.
\end{conjecture}
\end{svgraybox}
\noindent Let us say that a graph $G$ has the {\em AMRZ property\/} if $\lambda(G)+\mu(G)\ge \abs{G}$.
\begin{theorem}\label{thm:amrz+}{\rm\cite[corollary of Theorem 2.1]{AMRZ}}~
All graphs that can be constructed from chordal graphs 
by repeated substitutions and gluing vertices have the AMRZ property.
\end{theorem}
Theorem 2.1 of~\cite{AMRZ} is stronger than Theorem~\ref{thm:amrz+}: except for six graphs that have the AMRZ property, it replaces the AMRZ property by the property that the number of lines plus the number of bridges is at least the number of vertices.

In March 2018, Martin Matamala and Jos\' e Zamora proved that all bipartite graphs have the AMRZ property:
\begin{theorem}\label{thm:mzbip}{\rm\cite[corollary of Theorem 19]{MZ}}~
In all bipartite graphs except for $C_4$ and $K_{2,3}$, the number of lines plus the number of bridges is at least the number of vertices.
\end{theorem}
Theorem 19 of~\cite{MZ} is stronger than Theorem~\ref{thm:mzbip}: it replaces  the `number of lines' by `number of lines determined by pairs of vertices at distance 2'.
\begin{conjecture}\label{con:amrz+2}{\rm \cite[Conjecture 2.2]{AMRZ}}~
In all graphs with no pendant edges except for finitely many cases, the number of lines plus the number of bridges is at least the number of vertices.
\end{conjecture}
Here are Problems~\ref{prb:hhd}, \ref{prb:adj}, \ref{prb:non}, \ref{prb:dglue} with `DBE property' replaced by `AMRZ property' and phrased more cautiously:
\begin{svgraybox}
\begin{problem}\label{prb:hhd+} 
True or false? All HHD-free graphs have the AMRZ property.
\end{problem}
\end{svgraybox}\begin{svgraybox}
\begin{problem}\label{prb:adj+} 
True or false? Splitting a vertex into adjacent twins preserves the AMRZ property.
\end{problem}
\end{svgraybox}
\begin{svgraybox}
\begin{problem}\label{prb:non+} 
True or false? Splitting a vertex into nonadjacent twins preserves the AMRZ property.
\end{problem}
\end{svgraybox}
\begin{svgraybox}
\begin{problem}\label{prb:glu+} 
True or false? Gluing vertices preserves the AMRZ property.
\end{problem}
\end{svgraybox}
\noindent 

\section{Beyond graphs}\label{sec:bg}

\subsection{Metric spaces}\label{subsec:metr}
Just as all our graphs have at least two vertices, all our metric spaces have at least two points.

In the domain of metric spaces not necessarily arising from graphs, Conjecture~\ref{con:dbem} has been verified, in addition to its Euclidean case (Theorem~\ref{thm:predbe}), in another special case, that of nearly all finite subspaces of $({\bf R}^2,\ell_1)$. Here,`nearly all' means {\em non-degenerate\/} in the sense that no two points in the ground set share a coordinate.
\begin{theorem}\label{thm:kp1}{\rm \cite[Theorem 1.1]{KP13}}~
Every non-degenerate finite subspace of $({\bf R}^2, \ell_1)$ has the DBE property.
\end{theorem}
\begin{svgraybox}
\begin{problem}\label{prb:kp} 
Prove Theorem~\ref{thm:kp1} with the non-degeneracy assumption dropped.
\end{problem}
\end{svgraybox}
\begin{theorem}\label{thm:kp2}{\rm \cite[Theorem 1.2]{KP13}}~
Every finite subspace of $({\bf R}^2, \ell_1)$
has at least $\abs{V}/37$ distinct lines or a universal line. 
\end{theorem}
Problem~\ref{prb:kp} is a stepping stone toward proving that 
\begin{itemize}
\item[] (?) every finite subspace of every $({\bf R}^d, \ell_1)$ has the DBE property.
\end{itemize}
Allowing arbitrary values of $d$, but restricting the range of vectors in $V$ may seem to create another stepping stone, namely, proving that 
\begin{itemize}
\item[] (?) every finite subspace of every $(\{0,1\}^d, \ell_1)$ has the DBE property.
\end{itemize}
However, this restriction does not make the problem any easier: ($\{0,1\ldots ,k\}^d,\ell_1)$ is isometrically embeddable in $(\{0,1\}^{kd},\ell_1)$. To see this, allocate an ordered set of $k$ 
coordinates to each of the original $d$ coordinates and, within this set,
represent value $x$ by $1$s in the first $x$ positions followed by 0s in the
last $k-x$ positions.\\

There is nothing special about metric spaces with the $\ell_{\infty}$ metric: every metric space $(V, dist)$ is isometrically embeddable in $({\bf R}^{\abs{V}},\ell_{\infty})$. To see this, enumerate the elements of $V$ as 
$v_1,\ldots ,v_m$ and map each $v$ to $(dist(v,v_1),\ldots ,dist(v,v_m))$. (A related theorem of Fr\' echet \cite{Fre} states that every separable metric space is isometrically embeddable in the space of all bounded sequences of real numbers endowed with the supremum norm.) Since $({\bf R}^2,\ell_2)$ and $({\bf R}^2,\ell_{\infty})$ are isometric (one isometry maps $(x,y)$ to $(x+y,\,x-y)$), asking for a proof that all finite subspaces of $({\bf R}^2,\ell_{\infty})$ have the DBE property is just another way of stating Problem~\ref{prb:kp}.\\

Theorem~\ref{thm:glb} extends to metric spaces with a weaker lower bound:
\begin{theorem}\label{thm:mlb}{\rm \cite[Theorem 3.1]{ACH}}~
Every metric space on $n$ points has $\Omega(n^{1/2})$ distinct lines or a universal line.
\end{theorem}
Theorem~\ref{thm:d2} extends to metric spaces:
\begin{theorem}\label{thm:d2+}{\rm \cite[Theorem 1]{Chv14}}~
Every metric space on $n$ points  with distances in $\{0,1,2\}$ has the DBE property.
\end{theorem}
The notions of $\lambda$ and $\mu$ introduced in Subsection~\ref{subsec:inter} extend from graphs to metric spaces. In March 2018, Martin Matamala and Jos\' e Zamora~\cite{MZ} proved that the `DBE property' in Theorem~\ref{thm:d2+} can be replaced by `AMRZ property':
\begin{theorem}\label{tmm:d2++}{\rm\cite[Theorem 11]{MZ}}
Every metric space on $n$ points ($n\ge 3$) with distances in $\{0,1,2\}$ satisfies $\lambda+\max\{\mu -1, 0\}\ge n$.
\end{theorem}
A special case of Theorem~\ref{thm:asym3} extends to metric spaces:
\begin{theorem}\label{thm:masym}{\rm \cite[Theorem 5.3]{ACH}}~
Every metric space on $n$ points  with distances in $\{0,1,2,3\}$ has $\Omega(n^{4/3})$ distinct lines.
\end{theorem}
It is conceivable that the conclusion of this theorem remains valid even when the hypothesis is relaxed:
\begin{svgraybox}
\begin{problem}\label{prb:masym} Conjecture 1.3 of~\cite{ACH}:
Every metric space on $n$ points  with a constant number of distinct distances has $\Omega(n^{4/3})$ distinct lines.
\end{problem}
\end{svgraybox}
\noindent This conjecture, if valid, would subsume Theorem~\ref{thm:asym3} with constant $d(n)$. A partial result in its direction goes as follows:
\begin{theorem}\label{thm:diam+}{\rm\cite[Theorem 4.3]{ACH}}~
Every metric space on $n$ points ($n\ge 2$) with at most $k$ distinct nonzero distances has at least $n/5k$ distinct lines.
\end{theorem}

\subsection{Pseudometric betweenness}\label{subsec:pseudo}
To construct all lines in a prescribed metric space $M$, we need not know its distance function $dist$. The ternary relation $B(M)$ defined by  
\[(u,v,w)\in B(M) \;\Leftrightarrow\;
\mbox{ $u,v,w$ are all distinct and 
$ dist(u,v)+ dist(v,w)= dist(u,w)$ }
\]
suffices: lines in $M$ are determined by 
\begin{equation}\label{eq:line}
L(xy)\;=\;\{x,y\}\cup \{z:\; (x,y,z)\in B(M)\:\vee\: (y,z,x)\in B(M)\:\vee\: (z,x,y)\in B(M)\}.
\end{equation}
A ternary relation is called a {\em metric betweenness\/} if it is isomorphic to some $B(M)$.
Menger \cite{M28} seems to have been the first to study these 
relations. He pointed out that every metric betweenness $B$ has properties 
\begin{tabbing}
xx\=(M12)x\=\kill
\>(M0)\>$(u,v,w)\in B$ $\;\Rightarrow\;$ $u,v,w$ are three distinct points,\\
\>(M1)\>$(u,v,w)\in B$$\;\Rightarrow\;$ $(w,v,u)\in B$,\\ 
\>(M2)\>$(u,v,w)\in B$$\;\Rightarrow\;$ $(u,w,v)\not\in B$,\\
\>(M3)\>$(u,v,w),(u,w,x)\in B\;\Rightarrow\; (u,v,x),\;(v,w,x)\:\in B$.
\end{tabbing} 

Following~\cite{BBC13}, a ternary relation $B$ is called a {\em pseudometric betweenness\/} if it has properties (M0), (M1), (M2), (M3). Lines in a pseudometric betweenness $B$ are defined by by~(\ref{eq:line}) with $B$ in place of $B(M)$.

Euclidean betweenness is of course pseudometric and it has additional properties
\begin{tabbing}
xx\=(M12)x\=\kill
\>(M4)\>$(u,v,w),(v,w,x)\in B\;\Rightarrow\; (u,v,x),\;(u,w,x)\:\in B$,\\
\>(M5)\>$(u,v,w),(u,v,x)\in B\;\Rightarrow\; (u,w,x),(v,w,x)\in B \;\vee\; (u,x,w),(v,x,w)\in B$,\\
\>(M6)\>$(u,v,x),\;(u,w,x)\in B\;\Rightarrow\; (u,v,w),(v,w,x)\in B \;\vee\; (u,w,v),(w,v,x)\in B$,
\end{tabbing} 
but not every pseudometric betweenness with these properties is Euclidean: for instance,\\
\phantom{xxxxx}$\{(a_1,b_1,c_1),\;(a_1,b_2,c_2),\;(a_2,b_1,c_2),\;(a_2,b_2,c_1)$,\\ 
\phantom{xxxxx$\{$}$(c_1,b_1,a_1),\;(c_2,b_2,a_1),\;(c_2,b_1,a_2),\;(c_1,b_2,a_2)\}$\\
is not even metric. (For more on metric betweenness, see~\cite[Section 6]{Ch04}.)

It is conceivable that every  pseudometric betweenness has the DBE property. In February 2018, Pierre Aboulker [personal communication] proved this in the special case where the betweenness has a couple of additional properties and suggested that one of these two restrictions may be dropped:
\begin{theorem}\label{thm:pa}
Every  pseudometric betweenness with properties {\rm (M4)} and {\rm (M5)} has the DBE property.
\end{theorem}
\begin{svgraybox}
\begin{problem}\label{prb:pa}
True or false? Every  pseudometric betweenness with property (M4) has the DBE property.
\end{problem}
\end{svgraybox}
Theorem~\ref{thm:mlb} extends to pseudometric betweenness, although with an even weaker lower bound:
\begin{theorem}\label{thm:psd}{\rm \cite[Theorem 2.3]{ACH}}~
Every  pseudometric betweenness on $n$ points has $\Omega(n^{2/5})$ distinct lines or a universal line.
\end{theorem}

\subsection{$3$-uniform hypergraphs}\label{subsec:3uh}
To construct all lines in a prescribed pseudometric betweenness $B$, we need not know the order of the elements in each triple of $B$. The set  $T(B)$ of unordered triples defined by  
\[
T(B)=\{\{u,v,w\}:\,(u,v,w)\in B\}.
\]
suffices: lines in $B$ are determined by 
\begin{equation}\label{eq:hline}
L(xy)\;=\;\{x,y\}\cup \{z:\; \{x,y,z\}\in T(B)\}.
\end{equation}
Following~\cite{CheC}, {\em lines in a $3$-uniform hypergraph \/} with hyperedge set $T$ are defined by (\ref{eq:hline}) with $T$ in place of $T(B)$. Theorem~\ref{thm:predbe} cannot be generalized from the Euclidean plane all the way to $3$-uniform hypergraphs:
\begin{theorem}{\rm \cite[Theorem 3]{CheC}}~
There are arbitrarily large $3$-uniform hypergraphs with $n$ vertices , no universal line, and $\exp(O(\sqrt{\log n}))$ distinct lines.
\end{theorem}
Nevertheless, the number of distinct lines in $3$-uniform hypergraphs with $n$ vertices and no universal line grows beyond every bound as $n$ tends to infinity:
\begin{theorem}~{\rm \cite[Theorem 1]{ABC14}}
All $3$-uniform hypergraphs with $n$ vertices  have at least $(2-o(1))\lg n$ distinct lines or a universal line.
\end{theorem}
Here are four classes of $3$-uniform hypergraphs that are known to have the DBE property:
\begin{theorem}\label{thm:alm}{\rm \cite[Theorem 3]{ALM}}
If $H$ is a $3$-uniform hypergraph such that some graph $G$ shares its vertex set with $H$ and three vertices form a hyperedge in $H$ if and only if they are pairwise adjacent in $G$, then $H$ has the DBE property.
\end{theorem}
\begin{theorem}\label{thm.hyp}{\rm \cite[Theorems 2, 5, 6]{BBC13}}
  If, in a $3$-uniform hypergraph with at least two vertices,\\
\phantom{xxx} (a) no four vertices induce two hyperedges or\\
\phantom{xxx} (b) no four vertices induce one or three hyperedges or\\
\phantom{xxx} (c) no four vertices induce four hyperedges,\\
then the hypergraph has the DBE property.
\end{theorem}
By the way, Theorem 2 of ~\cite{BBC13} goes beyond the first part of Theorem~\ref{thm.hyp} by describing all $3$-uniform hypergraphs where no four vertices induce two hyperedges and the number of distinct lines equals the number of vertices. This theorem is a generalization of the real De Bruijn--Erd\H{o}s theorem (quoted here as Theorem~\ref{thm:dbe1}) since every family $\mathcal L$ of subsets of a set $V$ such that every two distinct points of $V$ belong to precisely one member of $\mathcal L$ is the family of lines in a $3$-uniform hypergraph where no four vertices induce two or three hyperedges. 

Theorem~\ref{thm.hyp} suggests the following questions:
\begin{svgraybox}
\begin{problem}\label{prb:ge2}~\cite[Question 2]{BBC13}
True or false? If, in a {\rm $3$}-uniform hypergraph, every
  sub-hypergraph induced by four vertices has at least two hyperedges,
  then the hypergraph has the DBE property.
\end{problem}
\end{svgraybox}

\begin{svgraybox}
\begin{problem}\label{prb:1or4}~\cite[Question 3]{BBC13}
True or false? If, in a {\rm $3$}-uniform hypergraph, every
  sub-hypergraph induced by four vertices has one or two or four
  hyperedges, then the hypergraph has the DBE property.
\end{problem}
\end{svgraybox}

Here is a counterpart of Problem~\ref{prb:geq} in the context of hypergraphs:
\begin{svgraybox}
\begin{problem}\label{prb:heq}
How difficult is it to recognize equivalence relations $\equiv$ such that 
$ab\equiv xy$ if and only if $L_H(ab)=L_H(xy)$ for some $3$-uniform hypergraph $H$?
\end{problem}
\end{svgraybox}

\noindent The following procedure rejects some of the candidates $\equiv$ partitioning edge sets of $K_n$ into classes $C_1$,\ldots $C_m$ (just like Algorithm G does), accepts some others (which Algorithm G never does), and gives up in the remaining cases.

\begin{tabbing}
{\bf Algorithm H:}\\
{\bf for}\, $i=1$ to $m$\, {\bf do}\, $L_i=$ the set of all endpoints of edges in $C_i$\, {\bf end}\\
{\bf while}\, \=there are pairwise distinct vertices $u,v,w$ and (not necessarily distinct)\\
\> subscripts $i,j$ such that $uv\in C_i$, $w\not\in L_i$, $vw\in C_j$, $u\in L_j$\\ 
{\bf do}\>add $w$ to $L_i$;\\ 
{\bf end}\\
{\bf if}\>$L_i\ne L_j$ whenever $i\ne j$\\
{\bf then}\> {\bf return}\, \= the hypergraph with hyperedge set consisting of all $\{u,v,w\}$\\ 
\>\>such that $u.v.w$ are pairwise distinct and $uv\in C_i$, $w\in L_i$ for some $i$;\\
{\bf else}\>{\bf if}\> there are distinct $i,j$ such that $\abs{L_i}=\abs{L_j}=n$\\
\>{\bf then}\>{\bf return} message {\sc Does not arise from any hypergraph};\\ 
\>{\bf else}\>{\bf return} message {\sc Don't know};\\ 
\>{\bf end}\\
{\bf end}
\end{tabbing}
For instance, given classes\ $C_1=\{14\}$, $C_2=\{24\}$, $C_3=\{34\}$, $C_4=\{12,23,13\}$, Algorithm H constructs $T=\{\{1,2,3\}\}$ and $L_1=\{1,4\}$, $L_2=\{2,4\}$, $L_3=\{3,4\}$, $L_4=\{1,2,3\}$, and so it returns the hypergraph with hyperedge set $T$. Given classes\ $C_1=\{15\}$, $C_2=\{25\}$, $C_3=\{35\}$, $C_4=\{45\}$, $C_5=\{12,23,34\}$, $C_6=\{13,24,14\}$, Algorithm H constructs $L_5=L_6=\{1,2,3,4\}$, \ldots , and so it returns message {\sc Don't know}.

Correctness of Algorithm H follows from the observation that, for all $3$-uniform hypergraphs $H$ such that $H(uv)=L_H(xy)$ if and only if $uv$ and $xy$ belong to the same $C_i$, its {\bf while} loop maintains the invariant $uv\in C_i \;\Rightarrow\; L_H(uv)\supseteq L_i$.

Of course, problems analogous to Problem~\ref{prb:heq} can be posed also with `metric spaces' or `pseudometric betweenness' in place of `$3$-uniform hypergraphs', but there we have nothing beyond Algorithm G.

\subsection{Recognition problems}\label{subsec:reco}

We have been discussing objects in a hierarchy of four levels:
\begin{enumerate}
\item metric spaces arising from graphs,
\item general metric spaces,
\item pseudometric betweenness,
\item $3$-uniform hypergraphs.
\end{enumerate}
Because every undirected graph $G$ gives rise to its metric space $M(G)$, every metric space $M$ gives rise to its pseudometric betweenness $B(M)$, and every pseudometric betweenness $B$ gives rise to its $3$-uniform hypergraph $H(B)$ with hyperedge set $T(B)$, this four-level hierarchy suggests six questions:
\begin{tabbing}
\phantom{xxl}\= Question 12: \= Does a prescribed metric space $M$ arise\\
\>\> from a graph $G$ as $M=M(G)$?\\
\> Question 13: \> Does a prescribed pseudometric betweenness $B$ arise\\
\>\> from a graph $G$ as $B=B(M(G))$?\\
\> Question 14: \> Does a prescribed $3$-uniform hypergraph $H$ arise\\
\>\> from a graph $G$ as $H=H(B(M(G)))$?\\
\> Question 23: \> Does a prescribed pseudometric betweenness $B$ arise\\
\>\> from a metric space $M$ as $B=B(M)$?\\
\> Question 24: \> Does a prescribed $3$-uniform hypergraph $H$ arise\\
\>\> from a metric space $M$ as $H=H(B(M))$?\\
\> Question 34: \> Does a prescribed $3$-uniform hypergraph $H$ arise\\
\>\> from a pseudometric betweenness $B$ as $H=H(B)$?
\end{tabbing}

\begin{itemize}
\item Question 12 is easy since only one graph $G$ may satisfy $M=M(G)$ with the prescribed metric space $M$: two vertices are adjacent in $G$ if and only if their distance in $M$ is $1$.
\item Question 13 is also easy since only one graph $G$ may satisfy $B=B(M(G))$ with the prescribed pseudometric betweenness $B$: vertices $u,w$ are adjacent in $G$ if and only if no $v$ satisfies $(u,v,w)\in B$.
\item Question 14: Let us call a $3$-uniform hypergraph $H$ {\em graphic\/} if there is a graph $G$ such that $H=H(B(M(G)))$. An induced sub-hypergraph of a graphic hypergraph may not be graphic. One example is the $H(B(M(G)))$ where $G$ consists of the cycle with edges $12$, $23$, $34$, $45$, $56$, $61$, and the additional vertex $7$ adjacent to the antipodal vertices $2$ and $5$. Here, the sub-hypergraph induced by the four vertices $1,2,3,5$ is not graphic.
\end{itemize}
\begin{svgraybox}
\begin{problem}\label{prb:14}
How difficult is it to recognize graphic hypergraphs?
\end{problem} 
\end{svgraybox}
\begin{itemize}
\item Question 23 can be answered in polynomial time: see~\cite[Section 6]{Ch04}.
\item Question 24: Following~\cite[Section 3]{BBC13}, let us call a $3$-uniform hypergraph {\em metric\/} if there is a metric space $M$ such that $H=H(B(M))$. Not all metric hypergraphs are graphic: one example is the $3$-uniform hypergraph with four vertices and one hyperedge. This hypergraph arises from the metric space with
\[
dist(a,b)=dist(b,c)=1,\;\;\; dist(a,c)=dist(a,d)=dist(b,d)=dist(c,d)=2
\]
and is not graphic. All induced sub-hypergraphs of a metric hypergraph are metric; a $3$-uniform hypergraph is called {\em minimal non-metric\/} if it is not metric, but all its proper induced sub-hypergraphs are. Three examples of minimal non-metric hypergraphs are given in ~\cite{BBC13}. Are there infinitely many minimal non-metric hypergraphs? 
\end{itemize}
\begin{svgraybox}
\begin{problem}\label{prb:24}
How difficult is it to recognize metric hypergraphs?
\end{problem}
\end{svgraybox}
\begin{itemize}
\item  Question 34: Following ~\cite[Section 3]{BBC13} again, let us call a  $3$-uniform hypergraph  {\em pseudometric\/} if there is a pseudometric betweenness $B$ such that $H=H(B)$. 
Not all pseudometric hypergraphs are metric: one example is the {\em Fano hypergraph,\/} whose hyperedges are the lines of the projective plane of order $2$. Like all $3$-uniform hypergraphs in which no two hyperedges share two vertices, it is pseudometric; since it does not have the Sylvester--Gallai property, it is not metric~\cite{Che}. All induced sub-hypergraphs of a pseudometric hypergraph are pseudometric; a $3$-uniform hypergraph is called {\em minimal non-pseudometric\/} if it is not pseudometric, but all its proper induced sub-hypergraphs are. The three examples of minimal non-metric hypergraphs given in ~\cite{BBC13} are also minimal non-pseudometric. Are there infinitely many minimal non-pseudometric hypergraphs?
\end{itemize}
\begin{svgraybox}
\begin{problem}\label{prb:34}
How difficult is it to recognize pseudometric hypergraphs?
\end{problem}
\end{svgraybox}
Here are counterparts of Problem~\ref{prb:LG} on the higher levels of the four-level hiearchy:

\begin{svgraybox}
\begin{problem}\label{prb:LM}
How difficult is it to recognize hypergraphs whose hyperedge set is the family of lines in some metric space?
\end{problem}
\end{svgraybox}
\begin{svgraybox}
\begin{problem}\label{prb:LB}
How difficult is it to recognize hypergraphs whose hyperedge set is the family of lines in some pseudometric betweenness?
\end{problem}
\end{svgraybox}
\begin{svgraybox}
\begin{problem}\label{prb:LH}
How difficult is it to recognize hypergraphs whose hyperedge set is the family of lines in some $3$-uniform hypergraph?
\end{problem}
\end{svgraybox}
Finally, one could also ask how difficult is it to recognize mappings $f:\binom{V}{2}\rightarrow 2^V$ such that some object on a specified level of the hierarchy has $L(uv)=f(uv)$ for all $uv$. However, this recognition problem is easy on the highest level: the following propositions are logically equivalent.
\begin{tabbing}
(A)\ \ \= Some $3$-uniform hypergraph $H$ has $L_H(uv)=f(uv)$ for all $uv$.\\
(B) \> If $u,v,w$ are pairwise distinct, then 
 $w\in f(uv) \;\Leftrightarrow\; v\in f(uw) \;\Leftrightarrow\; u\in f(vw)$.
\end{tabbing}
Furthermore, if (B) is satisfied, then the $H$ featured in (A) is unique (its hyperedges are all $\{u,v,w\}$ such that $u,v,w$ are pairwise distinct and $w\in f(uv)$), and so the recognition problem on each of the lower levels is Problem~\ref{prb:14} or Problem~\ref{prb:24} or Problem~\ref{prb:34}.

\section{Afterword}\label{sec: aw}

L\' aszl\' o Lov\' asz wrote ``It is easy to agree that if a conjecture is good, one expects that its resolution should advance our knowledge substantially.''~\cite{LL}. Would resolution of Conjecture~\ref{con:dbem} advance our knowledge substantially? No. Not unless you stretch the meaning of `substantially' enough to cover a theme that concerned fewer than two dozen people for the last dozen years. But there is something archetypal about the thrill of taking familiar concepts to unfamiliar territory. Think non-Euclidean geometry. Families of lines in graphs and metric spaces may never find applications comparable to those of hyperbolic geometry, but when we step through the looking glass to study them, what discoveries shall we make? Is Conjecture~\ref{con:dbeg} true or false? Conjecture~\ref{con:nonasym}? I want to know.

\bigskip

\begin{svgraybox}
Note added in proof: In March 2018, shortly after this paper was
submitted for publication, Xiaomin Chen and Ehsan Chiniforooshan
solved Problem 1: they proved  that all bisplit graphs with at least
80 vertices have the DBE property and verified the rest by computer
computations. In addition, they proved that all bisplit graphs with
$n$ vertices and no universal line have $\Omega(n^{4/3})$ distinct
lines. In August 2018, Laurent Beaudou, Giacomo Kahn, and Matthieu
Rosenfeld proved, without the use of a computer, that all bisplit
graphs have the DBE property:
\url{https://arxiv.org/pdf/1808.08710.pdf}
\end{svgraybox}

\vfill

\vfill

\noindent This is a correct \footnote{The project manager at SPi Global who handled the production of the book 
on behalf of Springer neglected to make four corrections requested by the author.} 
version of the paper published by Springer on pages 149--176 of the book 
{\em Graph Theory Favorite Conjectures and Open Problems - 2\/} edited by Ralucca Gera, Teresa W. Haynes, and Stephen T. Hedetniemi.\\
\url{https://www.springer.com/gp/book/9783319976846}

\end{document}